\documentclass{article}

\usepackage{xcolor}
\usepackage[a4paper, total={7in, 8in}]{geometry}
\usepackage[spanish, USenglish]{babel}
\usepackage[T1]{fontenc}
\usepackage{amsmath, amsthm, amsfonts, amssymb, stmaryrd}
\usepackage{mathrsfs}
\usepackage{graphicx}
\usepackage[affil-it]{authblk}
\usepackage{url}
\usepackage{subcaption} 
\usepackage{comment}
\usepackage{colortbl}
\usepackage{enumerate}
\usepackage{float} 
\usepackage[hidelinks]{hyperref} 

\newtheorem{thm}{Theorem}[section]

\newtheorem{cor}[thm]{Corollary}

\theoremstyle{definition}

\newtheorem{defn}[thm]{Definition}
\theoremstyle{remark}
\newtheorem{rem}[thm]{Remark}

\title{Degree–Based Topological Indices of a General Random Chain}
\author[1]{Saylé Sigarreta} 
\author[2]{Hugo Cruz-Suárez}
\author[3]{Sergio Torralbas Fitz}
\date{\today}

\affil[1,2]{Facultad de Ciencias Físico Matemáticas, Benemérita Universidad Autónoma de Puebla, Puebla, México}
\affil[3]{Department of Orthopedic Oncology, University of Miami Health System, Florida, United States}

\begin{document}
\maketitle

\begin{abstract}
In this paper, we examine a specific type of random chains and propose a unified approach to studying the degree-based topological indices, including their extreme values. We derive explicit analytical expressions for the expected values and variances of these indices and we establish the asymptotic behavior of the indices. Specifically, we analyze the first Zagreb index, Sombor index, Harmonic index, Geometric-Arithmetic index, Inverse Sum Index, and the second Zagreb index for various general random chains, including random phenylene, random polyphenyl, random hexagonal, and linear chains.
\end{abstract}

\textbf{MSC classes}: 05C50, 05C80

\section{Introduction}
A graph is defined as a pair $G=(V,E)$, where the elements of $V$ are the vertices of the graph $G$ and the elements of $E$ are its edges. The graphs considered in this manuscript are, unless otherwise specified, finite, simple, and connected. Undoubtedly, graph theory concepts are potentially applicable for many purposes. For instance, a chemical graph is a model of a chemical system, i.e., a graph is used to represent a molecule by considering the atoms as the nodes of the graph and the molecular bonds as the edges.

\medskip

\noindent
The topological indices, on the other hand, quantify the structural information contained in the graph and are independent of the numbering of the nodes and edges. In 1947, Wiener began to use topological indices to study some physico-chemical properties of alkanes, thus giving rise to chemical graph theory \cite{3}. In fact, since then, several topological indices have been introduced and extensively studied to better understand the molecular structure \cite{i1, E1, E2}. In particular, the Sombor index was recently introduced by Gutman in \cite{E3}. Recent papers about the Sombor index can be found in \cite{E4, E5, E6, r, r1, 0}. Nowadays, the characterization of molecular structures through topological indices remains a key focus of chemical graph theory, which plays a crucial role in designing molecules with specific physico-chemical or biological properties.

\medskip

\noindent
Many important topological indices, specifically degree-based indices, can be defined as

\begin{equation}\label{sa2}
   TI(G)= \sum_{vu\in E(G)} h(d_{v}, d_{u}),
\end{equation}
    
\noindent
where $h$ is some function with the property $h(x, y) = h(y, x)$ for $x, y \in \{1,2, \dots\}$ and $d_{v}$ is the degree of a node $v$. In the rest of the manuscript, we will denote $h(x,y)$ as $h_{x,y}$.

\noindent
In the field of random chains analysis, topological indices have been an evolving research topic for the past two decades \cite{shan}. Multiple topological indices have been analyzed for different random chains, such as random cyclooctane chains (\cite{c1}, \cite{y1}), random polyphenyl chains (\cite{p1}, \cite{p2}), random phenylene chains (\cite{ph1}, \cite{ph2}), random spiro chains (\cite{n6}, \cite{y2}) and random hexagonal chains (\cite{h2}, \cite{h1}). In the same vein, in \cite{25}, Definition 2.1, the concept of a general random chain was introduced. In this paper, we adopt this definition with slight modifications.

\begin{defn}\label{d1}
Given a graph \( G \), we say that a graph $H$ is a random chain generated by $G$ if the following conditions hold: $H$ contains  \( n \in \mathbb{N} \) copies of \( G \) (so, we can denote $H$ as $G_n$), any two copies of \( G \) are either adjacent (i.e., they are attached in a certain way) or disjoint (i.e., they have no common vertices) and the adjacent copies induce a path of $n$ vertices. The construction of such a random chain can proceed as follows:

    \begin{enumerate}[(a)]
        \item  $G_1 = G$ and $G_2$ consists of two copies of $G$ attached in a specific manner.
        \item For each $n > 2$, $G_n$ is constructed by attaching one copy of $G$ to $G_{n-1}$ in $m (m \geq 1)$ specific ways, resulting in the graphs $G_n^1, G_n^2, \cdots, G_n^m$ with probabilities $p_1, p_2, \cdots, p_m$ respectively, where $\sum_{i=1}^m p_i = 1$.
    \end{enumerate}
Hence, from now on, we will denote 
$G_n$ as $G(n, p_1, p_2, \ldots, p_m)$.

\end{defn}

\noindent
In addition, the article \cite{25} focused on establishing the distributions of Sombor indices in a general random chain, where explicit analytical expressions for the expected values and variances were derived. Note that, $G(n, p_1, p_2, \ldots, p_m)$ is constructed by a zero-order Markov process ($0MP$). It is important to highlight that random chains play crucial roles in chemistry and material science due to their unique structural and chemical properties. They are widely used in applications such as conducting polymers, organic electronics, elastomer design, and biocompatible materials (\cite{m1}, \cite{m2}).

\medskip

\noindent
Given $G(n, p_1, p_2, \ldots, p_m)$ as a random chain, let $L_i$ denote the link (attachment) selected at time $i \geq 2$, specifically, $L_2$ represents the initial link, and for $i \geq 3$, $L_i$ is a random variable with range $\{1,2, \dots ,m\}$ where $p_j = \mathbb{P}(L_i = j)$. For us, $G(n, p_1, p_2, \ldots, p_m)$ is referred to as a $0MP$-random chain with respect to a topological index $TI$, if for all $n \geq 3$
$$TI(G(n, p_1, p_2, \ldots, p_m)) - TI(G(n-1, p_1, p_2, \ldots, p_m)) = g(L_n),$$
where $g:\{1, 2, \dots ,m\} \to \mathbb{R}.$ In simple terms, the change in the calculation of the topological index from the time $n-1$ to $n$ is independent of the links selected in the previous steps.

\medskip

\noindent
The primary goal of this manuscript is to derive explicit formulas for the expected value, variance, and asymptotic distribution of a $0MP$-random chain with respect to the topological indices defined in Equation (\ref{sa2}). Furthermore, we show that several well-known random chains, including random phenylene, random polyphenyl, random cyclooctane, and linear chains, fit within this framework, which also allows us to use the same methodology to study the deterministic versions of these chains. As a result, various topological indices are examined for these structures, with several known results emerging as corollaries.

\section{Main Result}\label{s2}
\noindent
 In this section, we state and prove our main result. First, let $n\geq 2$ consider the following notation: $TI_{n}:=TI(G(n, p_1, p_2, \ldots, p_m))$, $\alpha:=\mathbb{E}(g(X))$ and $\beta:=\mathbb{V}(g(X))$, where $X \sim L_i$ for $i \geq 3.$ Here, when $m\geq 2$, we use $G\left(n, p_{1},  p_{2}, \dots, p_{m-1}\right)
 $ instead of $G\left(n, p_{1},  p_{2}, \dots, p_{m}\right)$ since $\sum_{i=1}^{m} p_i=1.$
\begin{thm}\label{tt1}
Let $G\left(n, p_{1},  p_{2}, \dots, p_{m-1}\right)$ be a $0MP$-random chain with respect to a topological index $TI$, then for $n \geq 2$ 
\begin{enumerate}[(a)]
    \item  $\mathbb{E}(TI_{n})=TI_{2}+\alpha(n-2).$
    \item $V(TI_{n})=\beta(n-2).$

    \item As $n \to \infty$, $\frac{TI_{n}-\alpha n }{\sqrt{\beta n}} \stackrel{D}{\longrightarrow}N(0,1)$.

 \item As $n \to \infty$, $\frac{TI_n}{n} \stackrel{a.s}{\longrightarrow} \alpha$.

\end{enumerate}

\end{thm}
\medskip

\begin{proof}Let $n \geq 3$, since $G(n, p_1, p_2, \ldots, p_m)$ is a $0MP$-random chain with respect to $TI$, it follows by definition that\\

\medskip
\begin{equation}\label{papa}
     TI_n-TI_{n-1}=g(L_n),
 \end{equation}
\medskip

\noindent
for some $g:\{1,2, \dots ,m\} \to \mathbb{R}.$ At this point,  by using the above recursive relation we have that
$$TI_n=TI_2+\sum_{i=3}^{n} g(L_i),$$
which directly allows us to obtain the expression  for the expectation and the variance. Likewise, since $\{g(L_i)\}_{i \geq 3}$ are independent and identically distributed (Chap 10, Theorem 10.4 \cite{y0}), due to the the Central Limit Theorem and the Strong Law of Large Number we can verified (c) and (d), respectively.
\end{proof}

\begin{rem}\label{oc}

At this point, it is important to emphasize that, given a $0MP$-random chain with respect to a topological index $TI$, the following equivalence holds: \( TI_{n} = TI_{2} + \alpha(n - 2) \) almost surely for \( n \geq 2 \) (a deterministic sequence) if and only if \( g(1) = g(2) = \dots = g(m) \); meaning that, the change is independent of the link selected at time \( n \) and remains constant throughout.  In particular, when $m=1$ by vacuity the previous equivalence is verified. Finally, note that, $TI_{n}=TI_{2}+a^{T} \cdot X$, where  $a^{T}=(g(1),g(2),\dots,g(m))$ and $X=(X_{1},X_{2},\dots,X_{m})$ is a multinomial random variable  with parameters $n-2$ and $(p_{1},p_{2},\dots,p_{m})$. In this context, it is useful to pointed out that, the approximation given in the previous theorem  for the convergence in distribution is identical to the one obtained by using the previous representation and the Central Limit Theorem in the case of random vectors \cite{27}. 

\end{rem}

\section{Application to 0MP-Random Chains}
The purpose of this section is to provide examples of $0MP$-random chains in relation to degree-based topological indices. Additionally, we will build on the results from the previous section, with the following corollaries encapsulating these concepts.

\begin{cor}\label{c2}
Let $RPC_{n}=RPC\left(n, p_{1}, p_{2}\right)$ be a random phenylene chain with $n \geq 2$ and $TI$ a degree-based topological index. Then
\begin{enumerate}[(a)]
\item $TI_{n}=AX+Bn+C$,

\item $\mathbb{E}\left(TI_{n}\right)=(Ap_{1}+B)n-2Ap_{1}+C$,

\item $\mathbb{V}\left(TI_{n}\right)=A^{2} p_{1}\left(1-p_{1}\right)(n-2)$,
\item As $n \to \infty$,
    $\frac{TI_{n}-\mathbb{E}\left(TI_{n}\right)}{\sqrt{\mathbb{V}\left(TI_{n}\right)}} \xrightarrow{D} N(0,1)$,

\item As $n \to \infty$,
  $\frac{TI_n}{n} \stackrel{a.s}{\longrightarrow} \mathbb{E}\left(TI_{n}\right)$,
\end{enumerate}
\hfill
\medskip

\noindent
where $A=2h_{2,3}-h_{2,2}-h_{3,3}$, $B=h_{2,2}+2h_{2,3}+5h_{3,3}$, $C=4h_{2,2}-6h_{3,3}$ and $X \sim Binomial(n-2,p_{1})$. In particular, given a phenylene chain \( rpc_n \) containing \( n \) hexagons, 
\[
TI(rpc_n) = A n_1 + Bn + C,
\]
where \( n_1 \) represents the number of type-1 links chosen up to the time \( n \).
\end{cor}

\begin{figure}
   \centering
        \includegraphics[width=0.5\textwidth]{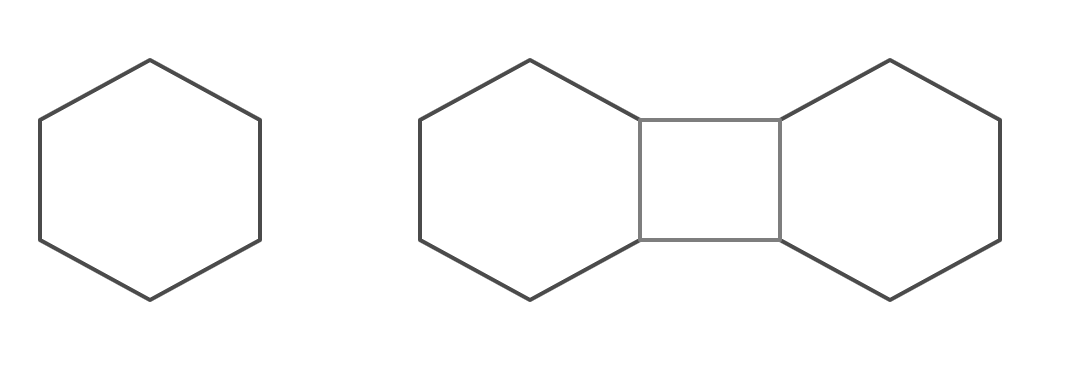}
    \caption{The graphs of $RPC_{1}$ and $RPC_{2}$.} 
    \label{f3}
\end{figure}

\begin{figure}
    \centering
    \includegraphics[width=1.2\textwidth]{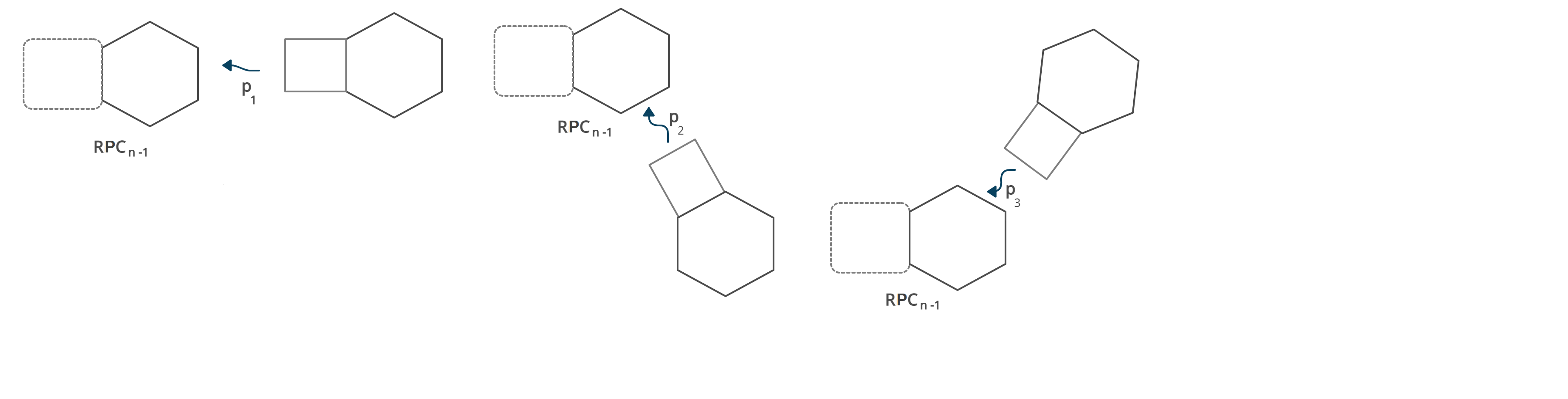}
   \caption{$RPC_{n}$.}
    \label{f4}
\end{figure}

\begin{proof}
A random phenylene chain is constructed by the following way: for $n=1$ and $n=2$, $RPC_{n}$ are shown in Figure \ref{f3}. For $n \geq 3$, the terminal hexagon can be attached in three ways, which results in $RPC_{n}^{1}, RPC_{n}^{2}$ and $RPC_{n}^{3}$, respectively, see Figure \ref{f4}. Given $RPC_{n}$, conducting the one-step analysis
 
\begin{eqnarray*}
RPC_{n}-RPC_{n-1}=\sum_{uv \in E_{n,1}} h_n(d_u,d_v) + \sum_{uv \in E_{n,2}} (h_n(d_u,d_v)-h_{n-1}(d_u,d_v)),
\end{eqnarray*}
\noindent  
where $  h_k(d_u,d_v):= h(d_u,d_v)$ is calculated within $RPC_{k}$,  $E_{n,1}$ represents the edges added as we progress from step $n-1$ to $n$ while $E_{n,2}$ are the edges that modify their $h(d_u,d_v)$ during the same transition. By analyzing each sum separately, we find that:

$$\sum_{uv \in E_{n,1}} h_n(d_u,d_v)=3h_{2,2}+2h_{2,3}+3h_{3,3},$$
and

\begin{align*}
\sum_{uv \in E_{n,2}} \left( h_n(d_u,d_v) - h_{n-1}(d_u,d_v) \right) 
&= h_{3,3} - 2h_{2,2} + h_{3,3}I_{\{L_n \neq 1\}} \notag \\
&\quad + \left( 2h_{2,3} - h_{2,2} \right)I_{\{L_n = 1\}}. 
\end{align*}

\noindent
Hence, by definition, $RPC_{n}$ is a $0MP$-random chain with respect to the degree-based topological indices with $g(1)=4h_{2,3}+4h_{3,3}$ and $g(2)=g(3)=h_{2,2}+2h_{2,3}+5h_{3,3}$. Therefore, by Theorem \ref{tt1}, (b)-(e) have been established. Finally, since $TI_{n}=TI_{2}+(g(1),g(2),g(2)) \cdot X$ for $X$ a multinomial random variable  with parameters $n-2$ and $(p_{1},p_{2},1-p_1-p_2)$, (a) is completed, and as a consequence, the deterministic formula is proven.
\end{proof}  

\begin{cor}\label{c3}
Let $RPoC_{n}=RPoC\left(n, p_{1}, p_{2}\right)$ be a random polyphenyl chain with $n \geq 2$ and $TI$ a degree-based topological index. Then
\begin{enumerate}[(a)]
\item $TI_{n}=AX+Bn+C$,

\item $\mathbb{E}\left(TI_{n}\right)=(Ap_{1}+B)n-2Ap_{1}+C$,

\item $\mathbb{V}\left(TI_{n}\right)=A^{2} p_{1}\left(1-p_{1}\right)(n-2)$,
\item As $n \to \infty$,
    $\frac{TI_{n}-\mathbb{E}\left(TI_{n}\right)}{\sqrt{\mathbb{V}\left(TI_{n}\right)}} \xrightarrow{D} N(0,1)$,

\item As $n \to \infty$,
  $\frac{TI_n}{n} \stackrel{a.s}{\longrightarrow} \mathbb{E}\left(TI_{n}\right)$,
\end{enumerate}
\hfill
\medskip

\noindent
where $A=h_{2,2}-2h_{2,3}+h_{3,3}$, $B=2h_{2,2}+4h_{2,3}+h_{3,3}$, $C=4h_{2,2}-4h_{2,3}-h_{3,3}$ and $X \sim Binomial(n-2,p_{1})$. Moreover, given a polyphenyl chain \( rpoc_n \) containing \( n \) hexagons, 
\[
TI(rpoc_n) = A n_1 + Bn + C,
\]
where \( n_1 \) represents the number of type-1 links chosen up to the time \( n \). 

\end{cor}

\begin{figure}
   \centering
    \includegraphics[width=0.5\textwidth]{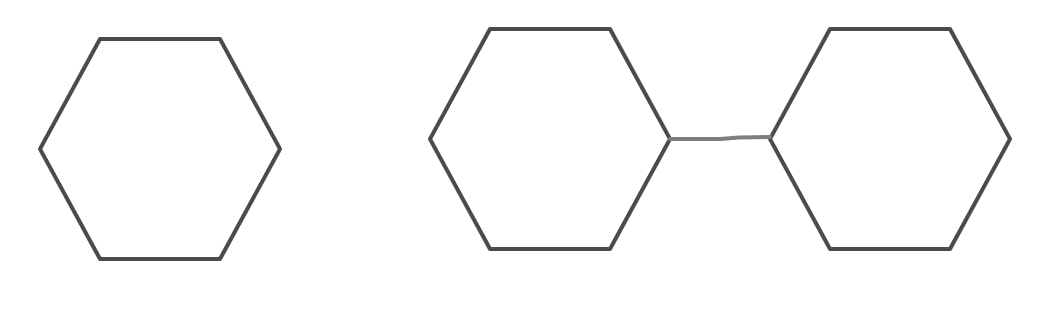}
    \caption{The graphs of $RPoC_{1}$ and $RPoC_{2}$.} 
    \label{f5}
\end{figure}

\begin{figure}
    \centering
    \includegraphics[width=0.95\textwidth]{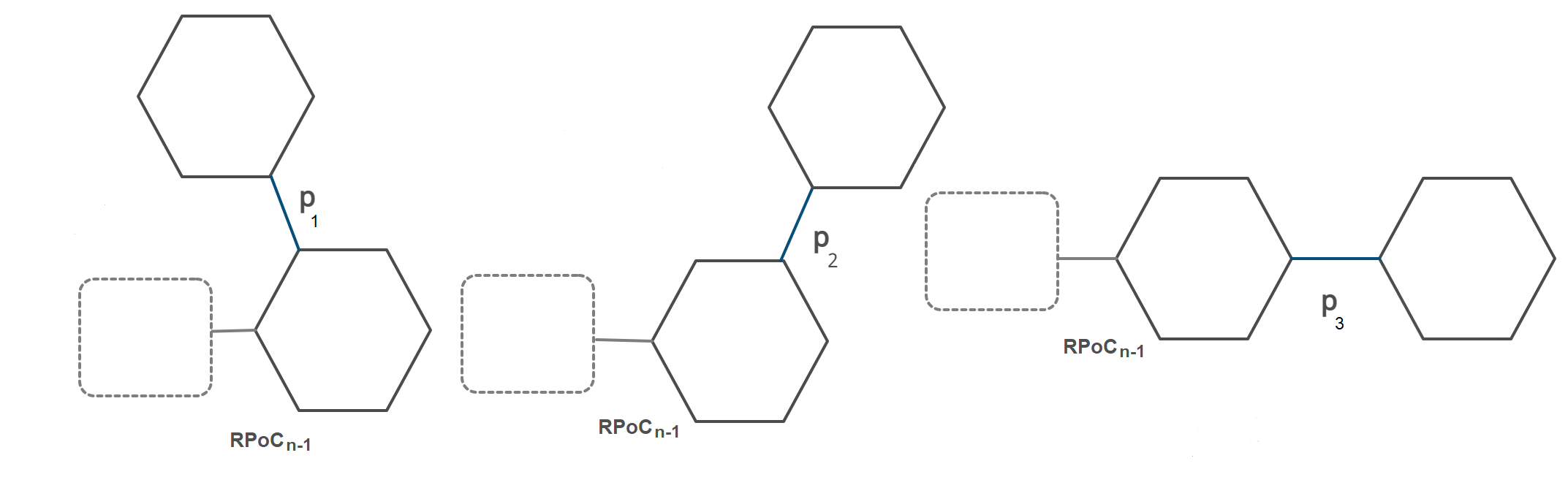}
   \caption{$RPoC_{n}$.}
    \label{f6}
\end{figure}

\begin{proof}
A random polyphenyl chain is constructed as follows: for $n=1$ and $n=2$, $RPoC_{n}$ are illustrated in Figure \ref{f5}. For $n \geq 3$, the terminal hexagon can be attached in three different ways, resulting in $RPoC_{n}^{1}, RPoC_{n}^{2}$ and $RPoC_{n}^{3}$, respectively, see Figure \ref{f6}. Following a similar method, we find that $RPoC_{n}$ is a $0MP$-random chain with respect to the degree-based topological indices with
$g(1)=2h_{2,2}+2h_{2,3}+2h_{3,3}$ and $g(2)=g(3)=2h_{2,2}+4h_{2,3}+h_{3,3}$. Therefore, applying the same procedure as in the proof of Corollary \ref{c2}, we complete the proof.
\end{proof}

\begin{cor}\label{c4}
Let $RHC_{n}=RHC\left(n, p_{1},p_{2}\right)$ be a random hexagonal chain with $n \geq 2$ and $TI$ a degree-based topological index. Then

\begin{enumerate}[(a)]
\item $TI_{n}=AX+Bn+C$,

\item $\mathbb{E}\left(TI_{n}\right)=(Ap_{1}+B)n-2Ap_{1}+C$,

\item $\mathbb{V}\left(TI_{n}\right)=A^{2} p_{1}\left(1-p_{1}\right)(n-2)$,
\item As $n \to \infty$,
    $\frac{TI_{n}-\mathbb{E}\left(TI_{n}\right)}{\sqrt{\mathbb{V}\left(TI_{n}\right)}} \xrightarrow{D} N(0,1)$,

\item As $n \to \infty$,
  $\frac{TI_n}{n} \stackrel{a.s}{\longrightarrow} \mathbb{E}\left(TI_{n}\right)$,

\end{enumerate}
\hfill
\medskip

\noindent
where $A=-h_{2,2}+2h_{2,3}-h_{3,3}$, $B=h_{2,2}+2h_{3,3}$, $C=4h_{2,2}+4h_{2,3}-2h_{3,3}$ and  $X \sim Binomial(n-2,p_{1})$. Moreover, given a hexagonal chain \( rhc_n \) containing \( n \) hexagons, 
\[
TI(rhc_n) = A n_1 + Bn + C,
\]
where \( n_1 \) represents the number of type-1 links chosen up to the time \( n \).

\end{cor}
\begin{figure}

   \centering
    \includegraphics[width=0.4\textwidth]{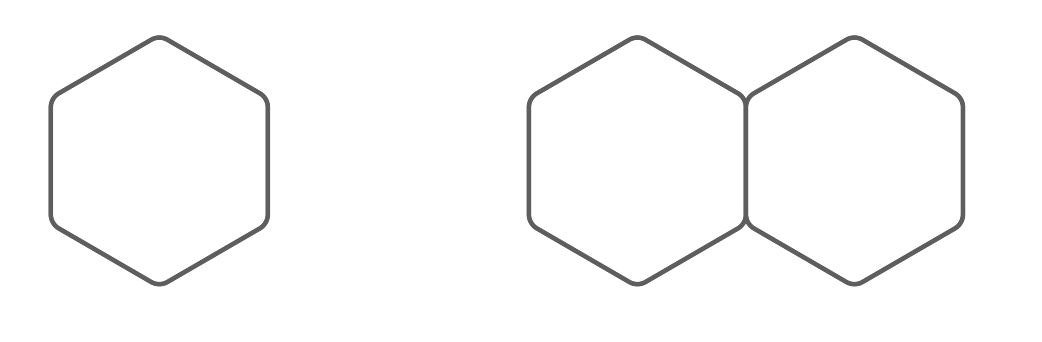}
    \caption{The graphs of $RHC_{1}$ and $RHC_{2}$.}

    \label{f7}
\end{figure}

\begin{figure}
    \centering
    \includegraphics[width=0.6\textwidth]{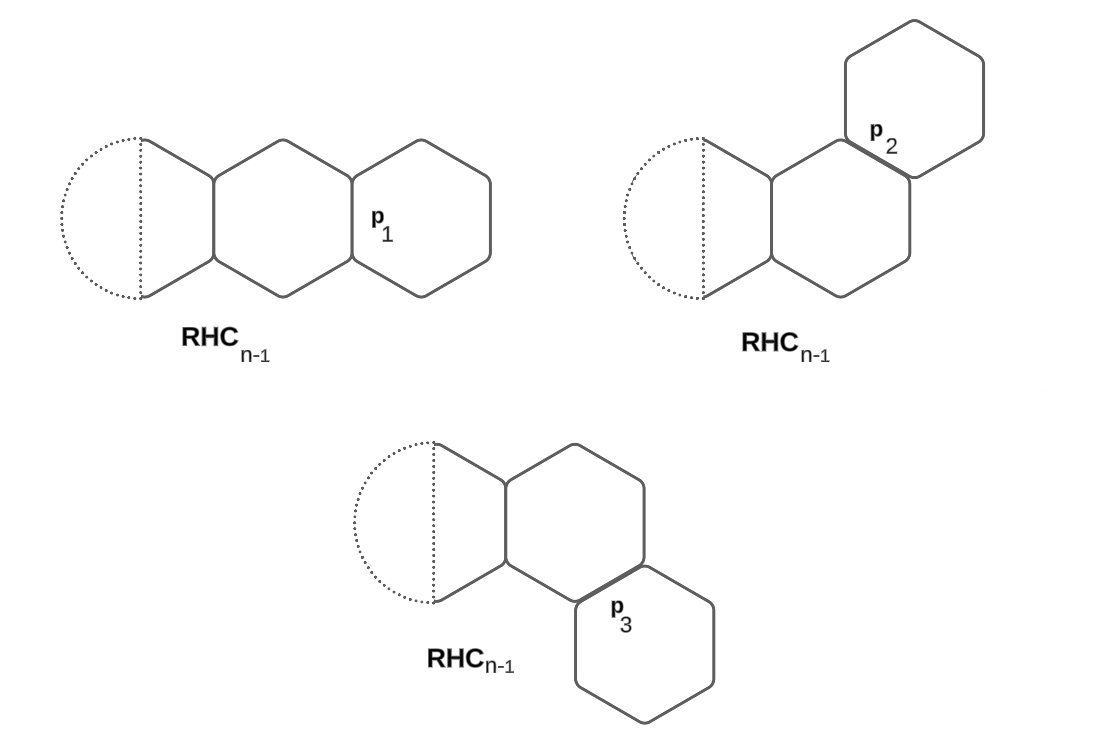}
   \caption{$RHC_{n}$.}
    \label{f8}
\end{figure}
\begin{proof}
A random hexagonal chain  is constructed by the following way: for $n=1$ and $n=2$, $RHC_{n}$ are shown in Figure \ref{f7}. For $n \geq 3$, the terminal octagon can be attached in three ways, which results in $RHC_{n}^{1}$ and $ RHC_{n}^{2}$ respectively, see Figure \ref{f8}. By following a similar method, we have that $RHC_{n}$ is a $0MP$-random chain with respect to the degree-based topological indices  with
$g(1)=2h_{2,3}+h_{3,3}$ and $g(2)=g(3)=h_{2,2}+2h_{3,3}$. Therefore, applying the same procedure of the proof of Corollary \ref{c2}, we complete the proof.

\end{proof}
\begin{rem}\label{r1}
By Remark \ref{oc}, we know that a degree-based topological index on the above random chains will be deterministic if and only if the corresponding A's are 0; which also is aligned with the expressions found in the previous corollaries. Likewise,  since $A_{RHC_n}=A_{R P o C_{n}}=-A_{R P C_{n}}$, it follows that: a fixed degree-based topological index is deterministic for $RHC_n$, $RPoC_n$ and $RPC_n$ if and only if
\begin{equation}\label{oc1}
    h_{2,2}+h_{3,3}=2 h_{2,3}.
\end{equation}
In particular, when $h(x, y)=x^{a}+y^{a}$ with $a \in \mathbb{R}$, i.e., we are working with the Generalized Zagreb Index, Equation (\ref{oc1}) holds. Consequently, we have $T I_{n}=B n+C=\displaystyle\sum_{v \in V(G)}\left(d_{v}\right)^{a+1}$, due to the identity \cite{S1} 
    \begin{center}
        $ \displaystyle\sum_{vu \in E(G)}\left(d_{v}\right)^{a}+\left(d_{u}\right)^{a}=
        \displaystyle\sum_{v \in V(G)}\left(d_{v}\right)^{a+1}$.
    \end{center} 
It is worth noting that  a similar argument may be used to demonstrate that random cyclooctane and random spiro chains \cite{y2} are also $0MP$-random chains with respect to the degree-based topological indices.

\end{rem}

\begin{rem} \label{r2}
Given the definition of a $0MP$-random chain within the context  of degree-based topological indices ($TI$), analyzing the extreme values of the function $g$, as determined by the recursive formula (\ref{papa}), could enables us to identify the extreme topological index values over the random chains. The analysis of the previously studied random chains ultimately reduces to a comparison between \( g(1) \) and \( g(2) \). If \( g(1) < g(2) \), then the minimum and maximum values of \( TI_n \) are achieved in \( G(n, 1, 0) \) and \( G(n, 0, p_2) \), respectively. Conversely, if \( g(1) > g(2) \), the minimum and maximum values of \( TI_n \) are attained in \( G(n, 0, p_2) \) and \( G(n, 1, 0) \), respectively.

\end{rem}

\begin{figure}[]

   \centering
    \includegraphics[width=0.5\textwidth]{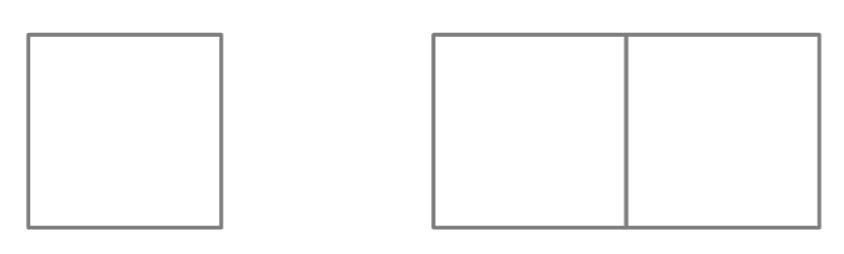}
    \caption{The graphs of $RSC_{1}$ and $RSC_{2}$.}

    \label{f9}
\end{figure}

\begin{figure}[h!]
    \centering
    \includegraphics[width=0.9\textwidth]{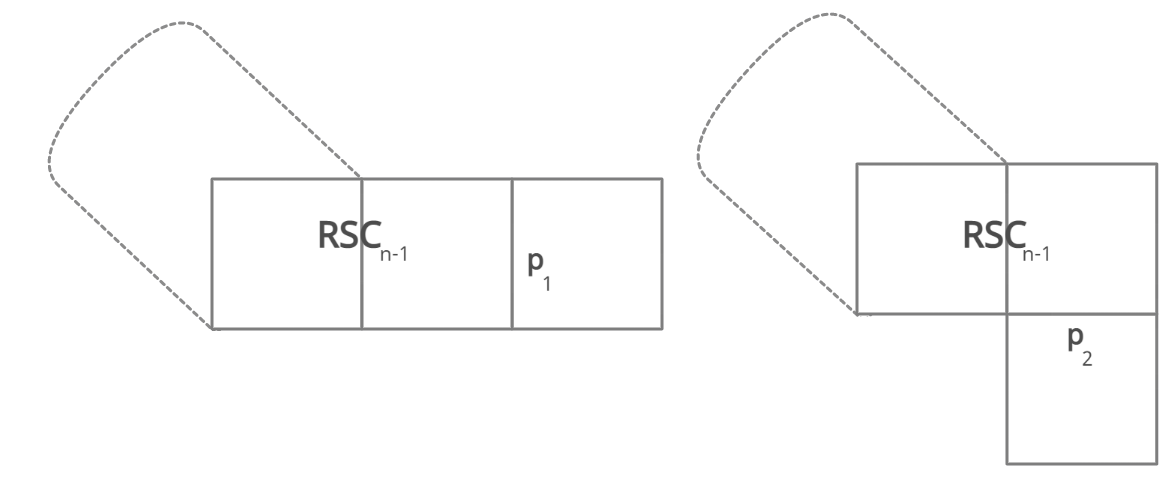}
   \caption{$RSC_{n}$.}
    \label{f10}
\end{figure}

\noindent
The authors in works such as (\cite{po5}, \cite{PP1}, \cite{PP2}, \cite{PP3}, \cite{PP4}) examined several well-known topological indices within the context of  a polyomino chain. A random polyomino chain at time $n \geq 1$ denoted as $RSC_{n}=RSC(n,p_{1},p_{2})$ can be constructed as follows: for $n=1$ and $n=2$, the configurations of $RSC_{n}$ are  illustrated in Figure \ref{f9}. For $n \geq 3$, the terminal square can be attached in two ways, resulting in $RSC_{n}^{1}$ and $RSC_{n}^{2}$ respectively, see Figure \ref{f10}. If $L_{n}=2$ for all $n \geq 3$ then $RSC_{n}=Z_{n}$ (a zigzag chain) whereas if  $L_{n}=1$ for all $n \geq 3$ then $RSC_{n}=Li_{n}$ (a linear chain). Note that, in general, \( RSC_{n} \) is not a \( 0MP \)-random chain with respect to degree-based topological indices. For example, the change between \( Z_{2} \) and \( Z_{3} \) differs from the change between \( Z_{4} \) and \( Z_{5} \) for the Randić index.
In contrast,  \(Li_n\) is classified as a \(0MP\)-random chain concerning degree-based topological indices, so, we can directly obtain the next corollary, which is aligned with Remark \ref{oc}. By the way, both \(Z_n\) and \(Li_n\) represent examples of general random chains with \(m=1\); however, while \(Li_n\) qualifies as a \(0MP\)-random chain regarding degree-based topological indices, \(Z_n\) does not. 

\begin{cor} \label{c5}
Let $Li_{n}$ be a linear chain with $n \geq 2$. Then

\begin{equation*}
    TI_{n}=An+B,
\end{equation*}

\hfill
\medskip

\noindent
where $A=3h_{3,3}$ and $B=2h_{2,2}+4h_{2,3}-5h_{3,3}$.

\end{cor}

Indeed, Corollary \ref{c5} is established in \cite{po5}. This analysis raises a natural question: Could a similar procedure be applied to specific types of topological indices in random polyomino chains? A preliminary exploration of this question is presented in \cite{y3,m3}.

\begin{figure}
   \centering
        \includegraphics[width=0.5\textwidth]{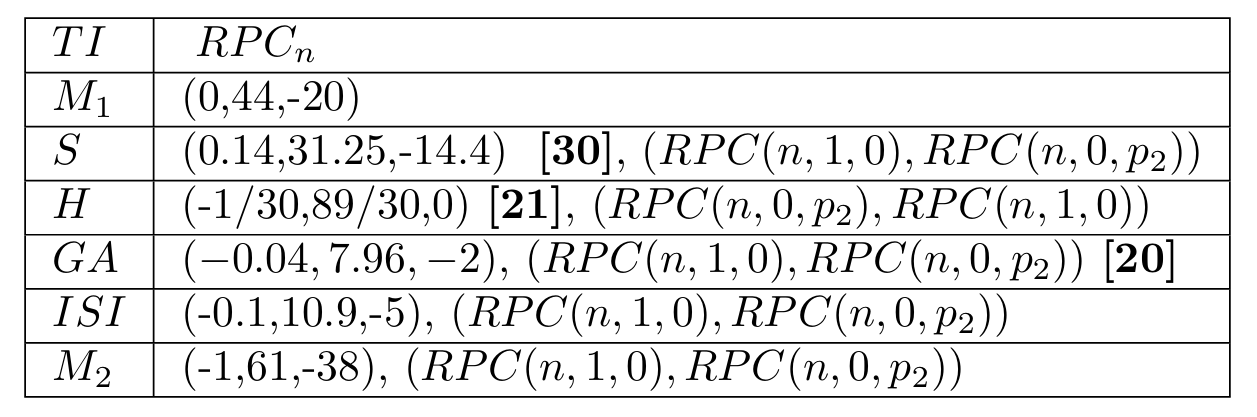}
    \caption{The relevant information associated with each topological index on $RPC_n$ is presented in the structure $\approx( A,B,C)$. The minimum and maximum values are respectively given in an ordered pair.} 
    \label{t1}
\end{figure}

\begin{figure}
   \centering
        \includegraphics[width=0.5\textwidth]{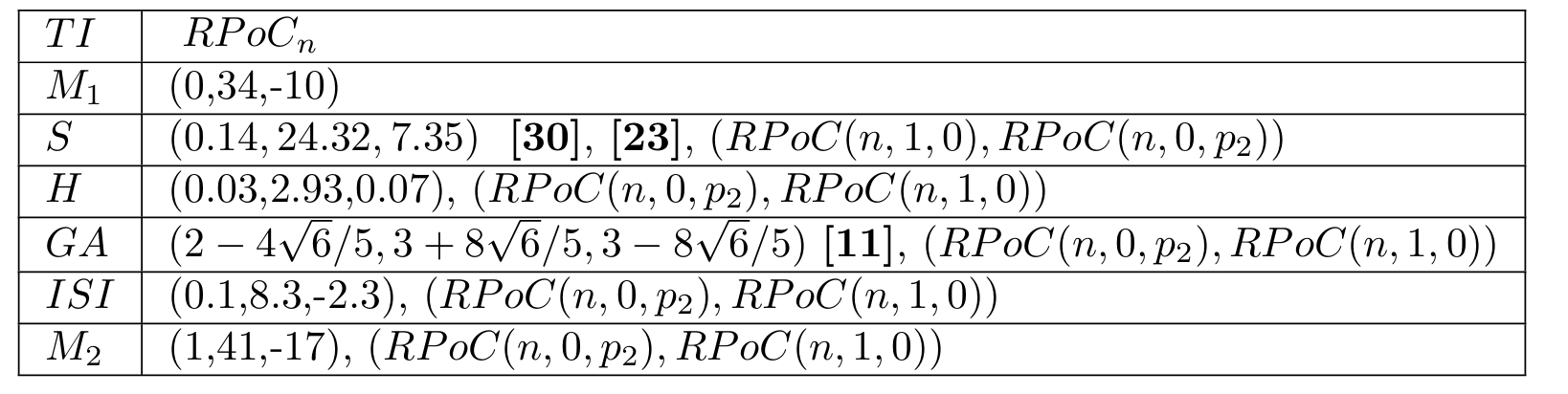}
    \caption{The relevant information associated with each topological index on $RPoC_n$ is presented in the structure $\approx( A,B,C)$. The minimum and maximum values are respectively given in an ordered pair.} 
    \label{t2}
\end{figure}

\begin{figure}
   \centering
        \includegraphics[width=0.5\textwidth]{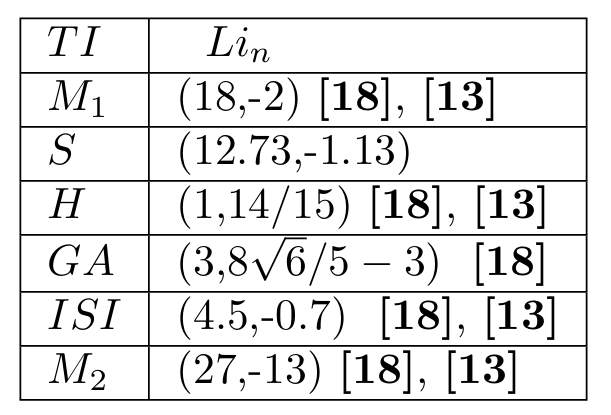}
    \caption{The information of interest associated with each topological index on $Li_n$ is exposed with the structure $\approx( A,B).$} 
    \label{t3}
\end{figure}

Finally, Figures \ref{t1}, \ref{t2}, \ref{t3}, and \ref{t3} display the values for the first Zagreb index ($M_{1}$), Sombor index ($S$), Harmonic index ($H$), geometric-arithmetic index ($GA$), inverse sum indeg index ($ISI$), and second Zagreb index ($M_{2}$) for the random chains analyzed in this section. Note that the results align with Remark \ref{r1} for $M_{1}$. Additionally, in light of Remark  \ref{r2}, these tables emphasize the chains where the extreme values of these topological indices are achieved.  In particular for these random chains the topological index of the chains at extreme values are attained can be directly computed as a consequence of the deterministic part of Corollaries \ref{c2}, \ref{c3}, and \ref{c4}. In fact, we have $TI(G(n,0,p_2)) = TI(G(n,0,1))$ since $g(2) = g(3)$.

By way of summary, in this paper,  we present a unified approach to studying degree-based topological indices in general random chains. We have derived the expected value, variance, and distribution of these indices. Additionally, we explored the asymptotic behavior and extreme values of the topological indices. Specifically, we focus on the Sombor, Harmonic, geometric-arithmetic, and inverse sum indeg indices for various general random chains, including random phenylene and random polyphenyl. Finally, we believe that the proposed method could be applied in future research to investigate additional characteristics, such as the average or deterministic expressions of both similar and different random chains, across a broader range of topological indices beyond those based solely on degree.

\begin{figure}
   \centering
        \includegraphics[width=0.5\textwidth]{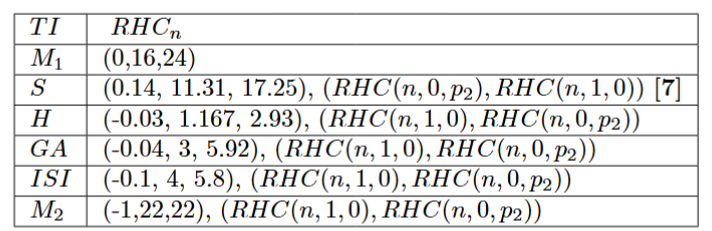}
    \caption{The relevant information associated with each topological index on $RHC_n$ is presented in the structure $\approx( A,B,C)$. The minimum and maximum values are respectively given in an ordered pair.} 
    \label{t4}
\end{figure}

\bibliography{Bib}
\bibliographystyle{acm}

\end{document}